\documentclass{amsart}

\usepackage[T1]{fontenc}
\usepackage{lmodern}

\usepackage{a4}

\usepackage[pdftitle={Freeness of equivariant cohomology and
    mutants of compactified representations},
  pdfauthor={Matthias Franz, Volker Puppe},
  pdfkeywords={torus actions, Hopf bundles,
    torsion-free equivariant cohomology}]{hyperref}

\let\epsilon\varepsilon

\def\N{{\mathbb N}}
\def\Z{{\mathbb Z}}
\def\Q{{\mathbb Q}}
\def\R{{\mathbb R}}
\def\C{{\mathbb C}}
\def\F{{\mathbb F}}
\def\m{\mathfrak m}
\def\L{\boldsymbol\Lambda}
\def\Lvee{\L^{\!\vee}}
\def\FF{{\mathcal F}}

\DeclareMathOperator{\rk}{rk}
\DeclareMathOperator{\coker}{coker}
\DeclareMathOperator{\Tor}{Tor}

\let\connsum\#
\let\connmult\star

\theoremstyle{plain}
\newtheorem{theorem}{Theorem}[section]
\newtheorem{proposition}[theorem]{Proposition}

\theoremstyle{definition}
\newtheorem{example}[theorem]{Example}
\newtheorem{remark}[theorem]{Remark}
\newtheorem{question}[theorem]{Question}

\newenvironment{acknowledgements}%
  {\begin{trivlist}\item[]{\em Acknowledgements.}\hskip0.5em}{\end{trivlist}}

\def\arxiv#1{preprint \href{http://arxiv.org/abs/#1}{\texttt{arXiv:#1}}}

\numberwithin{equation}{section}

\title[Freeness and mutants]%
  {Freeness of equivariant cohomology and\\
    mutants of compactified representations}

\author{Matthias Franz}
\address{Fachbereich Mathematik, 
  Universit\"at Konstanz, 78457 Konstanz, Germany}
\email{matthias.franz@ujf-grenoble.fr}

\author{Volker Puppe}
\address{Fachbereich Mathematik, 
  Universit\"at Konstanz, 78457 Konstanz, Germany}
\email{volker.puppe@uni-konstanz.de}



\begin{document}

\begin{abstract}
  We survey generalisations of the Chang--Skjelbred Lemma for
  integral coefficients. Moreover, we construct examples of
  compact manifolds with actions of tori of rank $> 2$
  whose equivariant cohomology is torsion-free, but not free.
  This answers a question of Allday's.
  The ``mutants'' we construct are obtained from compactified representations
  and involve Hopf bundles in a crucial way.
\end{abstract}

\maketitle

\section{Introduction}

Let $T=(S^1)^r$ be a torus and $X$ a ``sufficiently nice'' $T$-space,
for example, a compact (differentiable) $T$-manifold.
The equivariant cohomology~$H_T^*(X;\Q)$
is defined as the cohomology of
the Borel construction $X_T=ET\times_T X$. It captures quite a lot of
information about the $T$-action and sometimes also provides a link
to other subjects such as combinatorics, cf.~\cite{BuchstaberPanov:02}.
For that reason, one is interested
in efficient methods to compute $H_T^*(X;\Q)$.

The projection~$X_T\to ET/T=BT$ onto the classifying space of~$T$
gives $H_T^*(X;\Q)$ an algebra structure over
the polynomial ring~$R=H^*(BT;\Q)$. 
A very important special case is when $H_T^*(X;\Q)$ is a free $R$-module;
in this situation, one says that
``$X$ is (cohomologically) equivariantly formal'',
``$X_T$ has a cohomology extension of the fibre''
or just ``$X$ is CEF''.
Then, by a result of Chang--Skjelbred \cite[(2.3)]{ChangSkjelbred:74},
the sequence
\begin{equation}
  \label{C-S}
  0\longrightarrow
  H_T^*(X;\Q)\longrightarrow
  H_T^*(X^T;\Q)\stackrel\partial\longrightarrow
  H_T^{*+1}(X_1,X^T;\Q)
\end{equation}
is exact, where $X^T\subset X$ denotes the fixed point set,
$X_1$ the union of all orbits of dimension at most~$1$,
and $\partial$ the differential of the long exact cohomology sequence
for the pair~$(X_1,X^T)$.
In other words, $H_T^*(X;\Q)$ coincides, as subalgebra
of~$H_T^*(X^T;\Q)=H^*(X^T;\Q)\otimes R$,
with the image of~$H_T^*(X_1;\Q)\to H_T^*(X^T;\Q)$.
Usually $X^T$~and~$X_1$ are much simpler than~$X$,
so that \eqref{C-S} gives an easy method
to compute $H_T^*(X;\Q)$.
This can be used for instance for a short proof of
(a rational version of)
Jurkiewicz's description of
the cohomology of smooth projective toric varieties \cite{Jurkiewicz:80}
because they are known, like all Hamiltonian $T$-manifolds, to be CEF.

The purpose of this paper is twofold.
First we will survey generalisations of the Chang--Skjelbred sequence
for integer coefficients instead of the rationals,
as well as of a more general result due to Atiyah and Bredon.
This is motivated by the strong interest that
several participants of the Toric Topology Conference in Osaka
expressed for that topic. Secondly, we will answer
the following question raised by Chris Allday during his lecture
at the conference:

\begin{question}[Allday]
  Let $T=(S^1)^r$ be a torus of rank~$r>2$.
  Does there exist a compact $T$-manifold~$X$
  such that the $R$-module~$H_T^*(X;\Q)$ is torsion-free, but not free?
\end{question}

For rank~$r=1$ this is clearly false, even without assuming
Poincar\'e duality:
since $H_T^*(X;\Q)\subset H_T^*(X^T;\Q)=H^*(X^T;\Q)\otimes R$
is torsion-free and finitely generated,
it is free over the principal ideal domain~$R\cong\Q[t_1]$.
For~$r>1$ is it easy to find $T$-spaces~$X$
with $H_T^*(X;\Q)$ as required,
for example, the suspension~$\Sigma T$ of~$T$,
cf.~\cite[Ex.~5.5]{FranzPuppe:exact2}.
But $\Sigma T$ is not smooth for~$r>1$,
and when Allday posed the above question
(for rational Poincar\'e duality spaces),
he already proved that for~$r=2$
such spaces
cannot exist \cite[Prop.]{Allday:85}.

In this paper we will answer Allday's question in the affirmative
by exhibiting an example for~$r=3$.
Note that any example automatically gives further ones
for higher rank tori by adding circle factors
which act trivially on the space because in this case
\begin{equation}
  H_{T\times S^1}^*(X;\Q)=H_T^*(X;\Q)\otimes H^*(B S^1;\Q).
\end{equation}
Of course, this way one gets rather special modules
over polynomial rings.
Since our example turns out to be part of a (small) family of spaces
which seems to be of independent interest, we 
present the general construction.

Allday's question can be posed for
``$2$-tori''~$(\Z_2)^r$ as well, and,
in a spirit similar to~\cite{DavisJanuszkiewicz:91},
we will treat this case as well.

\begin{acknowledgements}
  One of us (M.\,F.) attended the Toric Topology Conference in Osaka
  and would like to thank the organisers for making
  such a pleasant and fruitful meeting possible, and also
  for financial support.
  Moreover, we are indebted to Jean-Claude Hausmann and Matthias Kreck
  for helpful discussions
  and to the referee for careful reading and useful suggestions
  which helped to improve the presentation.
\end{acknowledgements}

\section{Exact sequences for equivariantly formal
  \texorpdfstring{$T$}{T}-spaces}

Cohomology is taken with coefficients in~$\Z$
unless otherwise indicated. For any coefficients ring~$k$ (including~$\Z$),
we write $R=H^*(BT;k)=k[t_1,\dots,t_{\rk T}]$ with~$|t_i|=2$,
and $\m\subset H^*(BT;k)$ for the ideal of elements of positive degree.
All $R$-modules will be $\N$-graded.
Note that $R/\m=k$ is canonically an $R$-module
(concentrated in degree~$0$).
By the rank~$\rk M$ of an $R$-module~$M$ we mean
the dimension of the localisation of~$M$ over the quotient field of~$R$.
Tensor products without additional specification are taken over~$k$.

Let $X$ be a compact differentiable $T$-manifold or, more generally,
a finite $T$-CW~complex,
cf.~\cite[Def.~(1.1.1)]{AlldayPuppe:93}.
Denote by~$X_i$, $-1\le i\le r$, the equivariant $i$-skeleton of~$X$,
i.e., the union of all orbits of dimension~$\le i$.
In particular, $X_{-1}=\emptyset$, $X_0=X^T$ and $X_r=X$.
Each~$X_i$ is closed in~$X$.

The inclusion of pairs~$(X_i,X_{i-1})\hookrightarrow(X,X_{i-1})$
gives rise to a long exact sequence
\begin{equation}\label{exact-triple}
  \let\longrightarrow\to
  \longrightarrow H^*_T(X, X_i;k)
  \longrightarrow H^*_T(X, X_{i-1};k)
  \longrightarrow H^*_T(X_i, X_{i-1};k)
  \stackrel\delta\longrightarrow H^{*+1}_T(X, X_i;k)
  \longrightarrow,
\end{equation}
and likewise $(X_{i+1},X_{i-1})\hookrightarrow(X_{i+1},X_i)$ induces a map
\begin{equation}
  H^*_T(X_i, X_{i-1};k)\to H^{*+1}_T(X_{i+1}, X_i;k).
\end{equation}

Roughly at the same time as Chang and Skjelbred, Atiyah proved a much more
general theorem in the context of
equivariant $K$-theory~\cite[Ch.~7]{Atiyah:74}.
Bredon~\cite{Bredon:74} then observed that it applies equally to cohomology.
Bredon's version of Atiyah's result is the following:

\begin{theorem}[Atiyah--Bredon]\label{atiyah-bredon}
  If $H_T^*(X;\Q)$ is free over~$R$, then the sequence
  \begin{multline*} 
    \let\longrightarrow\to
    0
    \longrightarrow H^*_T(X;\Q)
    \longrightarrow H^*_T(X_0;\Q)
    \longrightarrow H^{*+1}_T(X_1, X_0;\Q)
    \longrightarrow \cdots
    \\
    \let\longrightarrow\to
    \cdots \longrightarrow H^{*+r-1}_T(X_{r-1}, X_{r-2};\Q)
    \longrightarrow H^{*+r}_T(X_r, X_{r-1};\Q)
    \longrightarrow 0
  \end{multline*}
  is exact.
\end{theorem}
A version for toric varieties
was proven
by Barthel--Brasselet--Fieseler--Kaup
for cohomology and intersection homology
\cite[Thm.~4.3]{BarthelBrasseletFieselerKaup:02},
this time including the (easier) converse.

In previous publications we obtained two variants
of Theorem~\ref{atiyah-bredon} for integer coefficients,
which we are going to recall now.
To state the first generalisation more succinctly,
we say that a closed subgroup~$K\subset T$ ``has at most one cyclic factor''
if the quotient~$K/K^0$ by the identity component is cyclic.
\begin{theorem}[\cite{FranzPuppe:exact1}]
  Assume that $H^*_T(X)$ is free over~$R$
  and that
  each isotropy group of~$X$ has at most one cyclic factor.
  Then the Atiyah--Bredon sequence is exact with integer coefficients.
\end{theorem}

We also obtained version
for other subrings of~$\Q$ and for prime fields.
The proof uses essentially the same techniques as in~\cite{Atiyah:74}.
Using the cohomological grading, which is absent in $K$-theory,
we could simplify the proof and obtain a  different version
which weakens the assumption on~$H_T^*(X)$ at the expense
of assuming more about the isotropy groups.
This time, we obtained an equivalence between various conditions.
(In fact, conditions \eqref{inclusion-fibre}~and~\eqref{Tor-1} below
are equivalent without any assumption on the isotropy groups.)

\begin{theorem}[\cite{FranzPuppe:exact2}]
  If all isotropy group of $X$ are connected,
  then the following conditions are equivalent:
  \def\theenumi{\roman{enumi}}
  \begin{enumerate}
  \item \label{inclusion-fibre}
    The inclusion of the fibre~$\iota\colon X\hookrightarrow X_T$
    induces a surjection~$\iota^*\colon H_T^*(X)\to H^*(X)$.
    Equivalently, the second map in the factorisation
    \begin{equation*}
      H_T^*(X)\to H_T^*(X)\otimes_R\Z\to H^*(X)
    \end{equation*}
    is an isomorphism.
  \item \label{Tor-1}
     $\Tor^{R}_1(H_T^*(X),\Z)=0$.
  \item \label{Atiyah-exact}
    The Atiyah--Bredon sequence is exact with integer coefficients.
  \end{enumerate}
\end{theorem}

In both cases, one can prove a variant of the Chang--Skjelbred Lemma
as well.

\begin{theorem}\label{C-S-Z}
  The Chang--Skjelbred sequence~\eqref{C-S} is exact over~$\Z$ if
  \def\theenumi{\roman{enumi}}
  \begin{enumerate}
  \item $H_T^*(X)$ is free over~$R$ and the isotropy group
    of each~$x\not\in X_1$ is contained in a proper subtorus, or
  \item $\Tor^{R}_1(H_T^*(X),\Z)=0$ and $T_x$ is connected for all~$x\in X_1$
    and contained in a subtorus of rank~$r-2$ for~$x\not\in X_1$.
  \end{enumerate}
\end{theorem}

This result is best possible as the examples
in \cite[Sec.~4]{TolmanWeitsman:99}~and~\cite[Sec.~5]{FranzPuppe:exact2} show.
In the setting of Hamiltonian group actions,
versions of the Chang--Skjelbred Lemma with integer coefficients
have been obtained by
Tolman--Weitsman~\cite[Prop.~7.2]{TolmanWeitsman:98},%
~\cite[Sec.~4]{TolmanWeitsman:99}
and Schmid~\cite[Thm.~3.2.1]{Schmid:01}
for connected as well as disconnected isotropy groups.
As mentioned in the introduction, the Chang--Skjelbred Lemma
(also in the version of Theorem~\ref{C-S-Z})
can be used to compute the equivariant cohomology of toric varieties
in cases where it is known to be free over~$R$.
This is explained in~\cite[Prop.~2.3]{BahriFranzRay:??}, for instance.

\section{An algebraic version of Allday's question}


If $X$ is a finite $T$-CW~complex such that $H^*(X;\Q)$
is a Poincar\'e duality algebra,
then there is a minimal Hirsch--Brown model
\begin{equation}\label{hirsch-brown}
  HB(X) = \bigl(H^*(X;\Q)\otimes R,\delta\bigr),
\end{equation}
see~\cite[Rem.~1.2.10, Rem.~3.5.9, Cor.~B.2.4]{AlldayPuppe:93}.
The cohomology of the differential $R$-module%
~\eqref{hirsch-brown} is the equivariant cohomology of~$X$.
The minimal Hirsch--Brown has the following properties:
\begin{enumerate}
\item\label{differential-linear}
  The image of the differential~$\delta$ lies in~$\m\cdot HB(X)$.
\item\label{bilinear-product}
  $HB(X)$
  carries an $R$-bilinear product~$\tilde\cup$
  (perhaps only associative and commutative up to homotopy).
  This product
  is compatible with~$\delta$ (meaning that $\delta$
  is a derivation of degree 1 with respect to it)
  and induces the cup product in cohomology.
  Moreover, it is a ``deformation'' of the
  cup product in~$H^*(X;\Q)$ in the sense that
  \begin{equation}
    H^*(X;\Q) = HB(X)\otimes_R\Q
  \end{equation}
  as $\Q$-algebras.
\item\label{duality-pairing}
  The composition of the product~$\tilde\cup$ in~$HB(X)$
  with the $R$-linear extension~$\tilde\sigma$
  of the orientation~$\sigma\colon H^*(X;\Q)\rightarrow\Q$
  gives a non-degenerate $R$-bilinear pairing
  \begin{equation}
    HB(X)\times HB(X)
    \stackrel{\tilde\cup}\longrightarrow HB(X)
    \stackrel{\tilde\sigma}\longrightarrow R,
  \end{equation}
  which is compatible with~$\delta$ and induces the
  Poincar\'e duality pairing on $H^*(X;\Q)$ upon tensoring with~$\Q$
  over~$R$.
\end{enumerate}

The above properties of the minimal Hirsch--Brown model
are essentially algebraic.
Thus, it is natural to consider
free differential $R$-modules
satisfying \eqref{differential-linear}--\eqref{duality-pairing}
and to ask whether in this context
torsion-free cohomology implies freeness.

Analysing Allday's proof for the case~$r=2$ 
shows that his result is purely algebraic in nature
and can essentially be stated in the following way:

\begin{proposition}[Allday]\label{allday-result}
  Let $\tilde C$ be a free differential 
  module 
  over $R=\Q[t_1,t_2]$
  with algebraic properties corresponding
  to (\ref{differential-linear})--(\ref{duality-pairing}) above.
  Then $H^*(\tilde C)$ is a free $R$-module
  if it is torsion-free.
\end{proposition}

On the other hand, as the following example shows, any finitely
generated $R$-module~$\tilde M$ can be realised as a
direct summand of~$H^*(\tilde C)$ for suitable~$\tilde C$.
Of course, it is not clear at
all whether this complex can be realised geometrically as the
minimal Hirsch--Brown model of some $T$-CW~complex.

\begin{example}\label{algebraic-example-1}
  Let $\tilde M$ be a finitely generated $R$-module, and
  choose a minimal free presentation
  \begin{equation}
    F_1\stackrel B\longrightarrow F_0\longrightarrow\tilde M.
  \end{equation}
  Define
  \begin{equation}
    \tilde C = R\oplus F_0\oplus F_1\oplus F_1'\oplus F_0'\oplus R[n]
  \end{equation}
  where ``$[n]$'' denotes a degree shift by~$n$,
  which is chosen large enough such that $F_i'\otimes_R \Q$
  is dual to~$F_i\otimes_R \Q$ for~$i=0$,~$1$.
  Choose a homogeneous basis of~$F_i$ and adjust the degrees of the
  elements in the dual basis of~$F_1'$ so that the
  degrees of an element and its dual add up to~$n$.
  The differential is given by~$B\colon F_1\to F_0$
  and its (graded) transpose~$B^T\colon F_0'\to F_1'$.
  The resulting cohomology is
  \begin{equation}
    H^*(\tilde C)
    = R\oplus\coker B\oplus\ker B\oplus\coker B^T\oplus\ker B^T\oplus R[n].
  \end{equation}
  In particular,
  $\tilde M = \coker B$ occurs as a direct summand
  in the cohomology of~$\tilde C$.
  The product on~$\tilde C$ is just the $R$-bilinear extension of the
  duality pairing on~$\tilde C\otimes_R \Q$.
  The differential defined above is compatible with this product.
\end{example}

The above result of Allday's implies in particular that
if for~$r = 2$ a non-free summand occurs in~$H^*(\tilde C)$,
then also torsion must occur.
The next example shows that Proposition~\ref{allday-result}
cannot be generalised to higher rank.

\begin{example}\label{algebraic-example-2}
  Let $r=3$ and set
  \begin{equation}
    \tilde C=R\oplus(R\oplus R\oplus R)[1]
    \oplus(R\oplus R\oplus R)[2]\oplus R[3].
  \end{equation}
  The differential~$\delta$ is zero
  except for a
  component~$\delta\colon(R\oplus R\oplus R)[2]\to(R\oplus R\oplus R)[1]$
  which is given by the matrix
  \begin{equation}
    B=\begin{bmatrix}
      0 & -t_3 & t_2 \\
      t_3 &    0 & -t_1 \\
      -t_2 & t_1 &    0
    \end{bmatrix}.
  \end{equation}

  The complex~$C=\tilde C\otimes_R \Q$ then is isomorphic to
  \begin{equation}
    \Q\oplus(\Q\oplus\Q\oplus\Q)[1]\oplus(\Q\oplus\Q\oplus\Q)[2]\oplus \Q[3]
  \end{equation}
  with trivial differential.
  The multiplication~$\tilde\mu\colon\tilde C\times\tilde C\to\tilde C$
  is the $R$-linear extension
  of the multiplication~$\mu\colon C\times C\to C$ 
  given by the standard dual pairings~$\Q\times \Q[3]\to \Q[3]$ and
  $(\Q\oplus \Q\oplus \Q)[1]\times(\Q\oplus \Q\oplus \Q)[2]\to \Q[3]$.

  In order to compute $H^*(\tilde C)$, we choose the canonical $R$-bases
  $(x_1,x_2,x_3)$ of $(\Q\oplus\Q\oplus\Q)[1]$ and
  $(y_1,y_2,y_3)$ of~$(\Q\oplus\Q\oplus\Q)[2]$.
  Then the kernel of~$B$ is generated by~$t_1 y_1+t_2 y_2+t_3 y_3$,
  which is of degree~$4$, and the cokernel is
  \begin{equation}\label{B-cokernel}
    (\Q\oplus\Q\oplus\Q)[1]\bigm/ \bigl\langle
      -t_3 x_2+t_2 x_3, t_3 x_1-t_1 x_3,-t_2 x_1+t_1 x_2
    \bigr\rangle.
  \end{equation}
  The assignment~$(\Q\oplus\Q\oplus\Q)[1]\to\Q[-1]$,
  $x_i\mapsto t_i$ induces an isomorphism of the quotient~\eqref{B-cokernel}
  with~$\m[-1]$. Hence,
  \begin{equation}
    H^*(\tilde C)=R\oplus\m[-1]\oplus R[3]\oplus R[4]
  \end{equation}
  as $R$-modules.
  In particular, $H^*(\tilde C)$ is torsion-free, but not free.
\end{example}

The above example cannot be
realised as the minimal Hirsch--Brown model of a finite
$T$-CW~complex~$X$: 
If this were the case, then, being torsion-free,
$H_T^*(X;\Q)$ would inject into~$H_T^*(X^T;\Q)=H^*(X^T;\Q)\otimes R$
by the Localisation Theorem, and the ranks over~$R$ would be equal.
More precisely, this would hold for even and odd degrees separately.
But since $\dim H_T^1(X;\Q)=3$ is greater than
\begin{equation}
  \dim H_T^1(X^T;\Q)\le \rk H_T^{\mathrm{odd}}(X^T;\Q)
    =\rk H_T^{\mathrm{odd}}(X;\Q)=2, 
\end{equation}
the restriction map cannot be injective in degree~$1$.

\begin{remark}\label{koszul-remark}
  The differential in the above example can be viewed as a part
  of the Koszul resolution of~$\Q$ over~$R$.
  (A good introduction to Koszul complexes can be found
  in~\cite[Sec.~17]{Eisenbud:95};
  see in particular Ex.~17.21 therein.)
  One can take a part of the Koszul resolution
  which is symmetric around the middle degree
  to obtain similar examples for $r>3$~variables.
  In the next section, we will construct geometric examples
  for certain values of~$r$, namely $r=3$,~$5$,~$9$,
  which realise these complexes (over~$\Z$ instead of~$\Q$) up to degree shift.
  Example~\ref{algebraic-example-2} is then realised
  by the minimal Hirsch--Brown model of the manifold~$Z_2$ for~$T=(S^1)^3$,
  see Section~\ref{S-1-equivariant}.
\end{remark}

\section{Mutants of compactified representations}\label{s-1-case}

Consider the standard action of $T=(S^1)^{r+1}$ on~$\C^{r+1}$,
which is given by
\begin{equation}
  (g_1,\dots,g_{r+1})\cdot(x_1,\dots,x_{r+1})=(g_1 x_1,\dots, g_{r+1}x_{r+1})
\end{equation}
for $g_i\in S^1$ and $x_i\in\C$. The action
extends to the one-point compactification~$S^{2r+2}$ of~$\C^{r+1}$
(which is a ``torus manifold'' in the sense of~\cite{HattoriMasuda:03}).
It has two fixed points, the origin~$0\in\C^{r+1}$
and the added point~$\infty$.
The union~$P=\R_{\ge0}^{r+1}\cup\{\infty\}$
of the positive orthant 
and the added point 
is topologically a cell~$D^{r+1}$ of dimension~$r+1$
and also a fundamental domain for the action.
In other words,
\begin{equation}
  S^{2r+2} \cong (D^{r+1}\times T)/\mathord\sim
\end{equation}
as topological spaces,
where $(x,g)\sim(x',g')$
if $x=x'$ and, in case $x\ne\infty$, $g_i=g'_i$ for each~$i$ with~$x_i>0$.

Now assume that~$r\in\{1,2,4,8\}$, and consider the Hopf bundle
\begin{equation}\label{hopf}
  S^{r-1}\hookrightarrow S^{2r-1}\stackrel p\longrightarrow S^r.
\end{equation}
Topologically, we define $Z_r$ to be the quotient
\begin{equation}\label{Z-r}
  Z_r=(D^{2r}\times T)/\mathord\sim,
\end{equation}
where again identification only takes place
at the boundary~$S^{2r-1}\subset D^{2r}$
and is induced by~$p$:
\begin{equation}
  (y,g)\sim(y',g')
  \;\Longleftrightarrow\;
  y,y'\in S^{2r-1} \;\text{and}\; \bigl((p(y),g)\sim(p(y'),g'\bigr).
\end{equation}

The sphere bundle
\begin{equation}\label{Y-r-X-r}
  S^{r-1}\hookrightarrow Y_r:=(\dot D^{2r}\times T)/\mathord\sim
  \longrightarrow X_r := (\dot D^{r+1}\times T)/\mathord\sim.
\end{equation}
is induced from the Hopf bundle~\eqref{hopf}
by a $T$-invariant map~$f\colon X_r\to X_r/T\cong\dot D^{r+1}\to S^r$.
In particular, it 
is orientable with $T$-invariant Euler class~$e$.
Note that the induced map of bundles is a retract,
so that we can consider $[S^r]\in H_r(X_r)$
and $[S^{2r-1}]\in H_{2r-1}(Y_r)$ in a canonical way.

The map~$f$ can be used to give $Z_r$
a smooth structure
in the following way:
The map~$\C^{r+1}\to\R_{\ge0}^{r+1}\subset\R^{r+1}$,
$(x_1,\ldots,x_{r+1})\mapsto(||x_1||^2,\ldots,||x_{r+1}||^2)$
is the quotient by~$T$ and extends
to~$S^{2r+2}\to P\subset S^{r+1}$.
Now choose a stereographic chart of~$S^{r+1}$ containing~$P$
and take the radial projection with respect to some inner point of~$P$.
The composition of all these maps can be used as~$f$,
and pulling back the smooth bundle~\eqref{hopf} along it
gives the smooth manifold~$Y_r$.
Because $Z_r$ is covered by $Y_r$
and~$((D^{2r}\setminus S^{2r-1})\times T)/\mathord\sim=(D^{2r}\setminus S^{2r-1})\times T$,
it is smooth, too.
Since $Z_r$ is in addition compact,
it satisfies Poincar\'e duality.

The $T$-action on~$Z_r$ is smooth
with fixed point set~$S^0\times S^{r-1}$
and quotient~$D^{2r}$.
Also note that for~$r\ne8$ the action of~$S^{r-1}$ on~$Y_r$
can be extended to~$Z_r$
by defining the complement~$Z_r\setminus Y_r\approx T$
to be fixed. The quotient of~$Z_r$ by this action is~$S^{2r+2}$.

\medbreak

We now calculate the integral homology of~$Z_r$.
In order to simplify this computation as well as that
for equivariant cohomology in Section~\ref{S-1-equivariant},
we will consider the action of~$\L=H_*(T)$ along the way.

Recall that $\L$ is an exterior algebra
with the Pontryagin product induced by the group multiplication.
It is generated by the classes~$x_i$
of loops around the different circle factors of~$T=(S^1)^{r+1}$.
Moreover, the action of $T$ on a space~$X$ induces
an action of~$\L$ on~$H_*(X)$ and also on~$H^*(X)$.
We will also need the quotient of~$\L$ by its top degree,
$\Lvee=\L/\L^{r+1}\cong\L^{<r+1}$.

Since $X_r$ is obtained from~$S^{2r+2}$ be removing one $T$-orbit
from the dense free stratum,
we have an exact sequence of $\L$-modules
\begin{equation}
  \cdots\longrightarrow H_{*+1}(S^{2r+2})
  \longrightarrow\L[r+1]_{*+1}
  \stackrel\partial\longrightarrow H_*(X_r)
  \longrightarrow H_*(S^{2r+2})
  \longrightarrow\cdots,
\end{equation}
hence an isomorphism of $\L$-modules
\begin{equation}\label{H-X-r}
  H_*(X_r)=\Z\oplus\Lvee[r],
\end{equation}
where 
the element~$1\in\Lvee[r]$ is mapped to~$[S^r]$.

Consider the Gysin homology sequence (cf.~\cite[Sec.~7]{Spanier:66})
of the sphere bundle~\eqref{Y-r-X-r},
\begin{equation}
  \cdots\longrightarrow H_{*+1-r}(X_r)
  \longrightarrow H_*(Y_r)
  \longrightarrow H_*(X_r)
  \stackrel{\cap e}\longrightarrow H_{*-r}(X_r)
  \longrightarrow\cdots,
\end{equation}
where $e$ is the Euler class. 
We claim that the map~$H_k(X_r)\to H_{k-r}(X_r)$
is an isomorphism for~$k=r$ and zero otherwise.
The first part follows by naturality
from the corresponding map for the Hopf bundle~\eqref{hopf}.
The second part uses that
capping with the $\L$-invariant class~$e$ is a $\L$-equivariant map
which sends the generator of~$\Lvee[r]$ to a $\L$-invariant.

Hence, we get a short exact sequence of $\L$-modules
\begin{equation}\label{sequence-lambda-Y}
  0\longrightarrow\Lvee[2r-1]
  \longrightarrow H_*(Y_r)
  \longrightarrow \Z\oplus\L^\diamond[r]
  \longrightarrow 0.
\end{equation}
Here $\L^\diamond$ denotes
the kernel of the map~$\Lvee\to\Z$ induced by the projection~$T\to 1$.
($\L^\diamond$ is the direct sum of~$\L^k$ for~$0<k<r+1$.)
The element~$1\in\Lvee[2r-1]$ is mapped to~$[S^{2r-1}]$.
The sequence~\eqref{sequence-lambda-Y} actually splits,
but it requires some work to see this, cf.~Remark~\ref{quasi-iso} below.

Since $Z_r$ is obtained from~$Y_r$ by gluing in an equivariant
$T$-cell of dimension~$2r$ along~$\dot D^{2r}\times T$,
we get an exact sequence of $\L$-modules
\begin{equation}
  \cdots\longrightarrow\L[2r]_{*+1}
  \stackrel\partial\longrightarrow H_*(Y_r)
  \longrightarrow H_*(Z_r)
  \longrightarrow\L[2r]_*
  \longrightarrow\cdots
\end{equation}
where the element~$1\in\L[2r]$ is also mapped to~$[S^{2r-1}]$,
which generates $\Lvee[2r-1]\subset H_*(Y_r)$.
By $\L$-equivariance, the sequence therefore splits into
the exact sequence
\begin{equation}\label{sequence-lambda-Z}
  0\longrightarrow\Lvee[2r-1]
  \longrightarrow H_*(Y_r)
  \longrightarrow H_*(Z_r)
  \longrightarrow\Z[3r+1] 
  \longrightarrow 0.
\end{equation}
Because the same submodule of~$H_*(Y_r)$ appears
in both \eqref{sequence-lambda-Y}~and~\eqref{sequence-lambda-Z}
and moreover $\Z[3r+1]$ is the only contribution in degrees~$\ge 3r+1$,
we finally obtain
an isomorphism of $\L$-modules
\begin{equation}\label{H-Z-r}
  H_*(Z_r)=\Z\oplus\L^\diamond[r]\oplus\Z[3r+1].
\end{equation}

Note that
$Z_1$ has the homology of~$S^2\times S^2$, and
for~$r\in\{2,4,8\}$ the homology of~$Z_r$ is (additively)
that of the connected sum
\begin{equation}
  \binom{r+1}1\connmult\bigl(S^{r+1}\times S^{2r}\bigr)
  \:\connsum\:\cdots\:\connsum\:
  \binom{r+1}{r/2}\connmult\bigl(S^{3r/2}\times S^{3r/2+1}\bigr),
\end{equation}
where
``\,$n\connmult X$\,'' means taking $n$~copies of the space~$X$
in the sum.
In Section~\ref{identify-spaces} we will identify $Z_1$~and~$Z_2$.

\section{Computing the equivariant cohomology}\label{S-1-equivariant}

By replacing each $T$-space~$X$ by its Borel construction~$ET\times_T X$,
one could compute the equivariant cohomology of~$Z_r$
in a way analogous to the calculation of the homology
in Section~\ref{s-1-case}.
Because this time the extension problems one is faced with
are more intricate, we follow a different approach which makes use of
the $\L$-structure on~$H^*(Z_r)$.

Consider the Leray--Serre spectral sequence
of principal $T$-bundle~$ET\times_T Z_r\to BT$.
Its $E_2$-term is of the form
\begin{equation}\label{E-2-S-1}
  E_2^{p,q}=R^p\otimes H^q(Z_r),
  \quad
  d_2(s\otimes\gamma)=\sum_{i=1}^{r+1} t_i s\otimes x_i\cdot\gamma,
\end{equation}
see for example \cite[Sec.~5.1]{Franz:03}.

Given the form~\eqref{H-Z-r} of~$H^*(Z_r)$,
\eqref{E-2-S-1} is essentially the Koszul resolution
of~$\Z$ over~$R$, with lowest and highest degree moved apart
from the central piece as in Remark~\ref{koszul-remark}.
We therefore obtain an isomorphism of $R$-modules
\begin{equation}
  E_3 =
  \begin{cases}
    R\oplus R[2]\oplus R[2]\oplus R[4] & \text{if $r=1$,}\\
    R\oplus \m[r-1]\oplus R[2r+2]\oplus R[3r+1] & \text{if $r\in\{2,4,8\}$.}
  \end{cases}
\end{equation}
In all cases, the rank of the $E_3$-term over~$R$ is $4$.
Since $E_3$ is torsion-free, any higher differential would lower the rank.
By the Localisation Theorem, the rank of~$H_T^*(Z_r)$ is the
same as that of
\begin{equation}
  H_T^*(Z_r^T)=H^*(S^0\times S^{r-1})\otimes R
\end{equation}
which is again $4$.
Hence, higher differentials cannot not occur.

Clearly, there is no extension problem for~$r=1$.
For~$r>1$, there is none either because of the equality
\begin{equation}
  H_T^k(Z_r)=\m[r-1]_k
\end{equation}
for~$k=r+1$~and~$r+3$:
the generators~$m_i$ of~$\m[r-1]$ live in degree~$r+1$,
and the relations~$t_i m_j=t_j m_i$ for~$i\ne j$ in degree~$r+3$.
We therefore get an isomorphism of $R$-modules
\begin{equation}
  H_T^*(Z_r) \cong
  \begin{cases}
    R\oplus R[2]\oplus R[2]\oplus R[4] & \text{if $r=1$,}\\
    R\oplus \m[r-1]\oplus R[2r+2]\oplus R[3r+1] & \text{if $r\in\{2,4,8\}$.}
  \end{cases}
\end{equation}
In particular, $H_T^*(Z_r)$ is free over~$R$ for~$r=1$
and torsion-free, but not free for larger~$r$.

\begin{remark}\label{quasi-iso}
  Alternatively, one can
  prove that the isomorphism~\eqref{H-Z-r} is induced by a
  quasi-isomorphism of $\L$-modules
  \begin{equation}
    \Z\oplus\L^\diamond[r]\oplus\Z[3r+1]\to C_*(Z_r)
  \end{equation}
  where $C_*(\cdot)$ denotes the normalised singular chain functor.
  An element $a\in\L^\diamond[r]$ is mapped to~$a\cdot c$,
  where $c\in C_r(Y_r)\subset C_r(Z_r)$ is a suitable
  transgression chain of the fibration~$Y_r\to X_r$
  and the action of~$\L$ is lifted from~$H_*(Z_r)$ to~$C_*(Z_r)$.
  Using the ``singular Cartan model'' (see~\cite[Sec.~5.1]{Franz:03}),
  this implies that the
  differential $R$-module~\eqref{E-2-S-1} computes $H_T^*(Z_r)$
  without extension problem.
\end{remark}

\section{The analogous construction for
  \texorpdfstring{$2$}{2}-tori}\label{z-2-case}

By replacing $S^1$ by~$\Z_2$ in the previous constructions,
one arrives at a similar family of spaces, which we 
describe in more detail in this section.
We assume again $r\in\{1,2,4,8\}$ and
consider 
the canonical action of the ``$2$-torus''~$G=(\Z_2)^{r+1}$
on the one-point compactification~$S^{r+1}$ of~$\R^{r+1}$
(which may be called a ``$2$-torus manifold'').
As before, a fundamental domain is the compactified positive orthant,
\begin{equation}
  S^{r+1}\cong(D^{r+1}\times G)/\mathord\sim,
\end{equation}
where identification only takes place along the boundary of~$D^{r+1}$.

Using the Hopf bundle~$p$ \eqref{hopf}, we define
\begin{equation}
  Z_r=(D^{2r}\times G)/\mathord\sim
\end{equation}
with the identification induced by~$p$.
Also define $Y_r$ and $X_r$ analogously to
the constructions in
Section~\ref{s-1-case}.
$Z_r$ is smooth by a reasoning similar to the previous one,
and there is again a $G$-invariant retract from $Y_r\to X_r$
to the Hopf bundle.
The $G$-action on~$Z_r$ has
fixed point set~$S^0\times S^{r-1}$ and quotient $D^{2r}$.
For~$r\in\{1,2,4\}$
the quotient of~$Z_r$ by the action of~$S^{r-1}$
is $S^{r+1}$.

To describe the homology of~$Z_r$, we consider
the group algebra~$A=H_*(G)=\Z[G]$ and its quotient~$A^\vee$
by the ``top element''
\begin{equation}
  \omega = (1-g_1)\cdots(1-g_{r+1}),
\end{equation}
where the~$g_i$ are the canonical generators of~$G$.
(Note that the line through $\omega$ is $A$-stable.)
Also let $A^\diamond$ be the kernel
of the augmentation~$A\to1$, divided by~$\Z\,\omega$.

Since $\omega\cdot[S^r]=0\in H_r(X_r)$, we get
\begin{equation}
    H_*(X_r) = \Z\oplus A^\vee[r],
\end{equation}
where $1\in A^\vee[r]$ corresponds to~$[S^r]$.
The Gysin homology sequence splits again into a short exact sequence,
\begin{equation}
   0\longrightarrow A^\vee[2r-1]
  \longrightarrow H_*(Y_r)
  \longrightarrow \Z\oplus A^\diamond[r]
  \longrightarrow 0
\end{equation}
because the Euler class~$e$ of the bundle is killed by all elements~$1-g_i$.

Arguing as before, we get an exact sequence of $A$-modules,
\begin{equation}
  0\longrightarrow A^\vee[2r-1]
  \longrightarrow H_*(Y_r)
  \longrightarrow H_*(Z_r)
  \longrightarrow\Z\,\omega[2r] 
  \longrightarrow 0,
\end{equation}
and finally
\begin{equation}
   H_*(Z_r)=\Z\oplus A^\diamond[r]\oplus\Z\,\omega[2r].
\end{equation}
Note that $Z_r$ has the homology of a connected sum of $2^r-1$~copies
of~$S^r\times S^r$.

\medbreak

The group algebra~$H_*(G;\F_2)$
is a strictly exterior algebra on generators~$1-g_i$ of degree~$0$,
and $R=H^*(BG;\F_2)$ a polynomial algebra on generators~$t_i$
of degree~$1$.
Also observe that all~$g_i$ act trivially on~$\F_2\,\omega$.

We claim that the minimal Hirsch-Brown model of the $G$-space~$Z_r$
is given by
\begin{equation}
  H^*(Z_r;\F_2)\otimes R,
  \quad
  \delta(\gamma\otimes s)=\sum_{i=1}^{r+1} (1-g_i)\cdot\gamma\otimes t_i s.
\end{equation}
Since $Z_r^G$ is not empty, the differential~$\delta$ does not hit
$H^{0}(Z_r;\F_2)\otimes R$,
and $\delta$ vanishes on $H^{2r}(Z_r;\F_2)\otimes R$
because of Poincar\'e duality,
see~\cite[Cor.~(5.3.4)]{AlldayPuppe:93}. So, for degree reasons,
$\delta$ only contains terms linear in the~$t_i$'s,
and these are given by the induced action in cohomology,
cf.~\cite[p.~453]{AlldayPuppe:93}. This proves the claim.

We therefore get
\begin{align}
  H_G^*(Z_r;\F_2) \cong
  \begin{cases}
    R\oplus R[1]\oplus R[1]\oplus R[2] & \text{if $r=1$,}\\
    R\oplus \m[r-1]\oplus R[r+1]\oplus R[2r] & \text{if $r\in\{2,4,8\}$.}
  \end{cases}
\end{align}
Again, $H_G^*(Z_r;\F_2)$ is free over~$R$ for~$r=1$
and torsion-free, but not free for larger~$r$.

\section{Identifying some mutants}\label{identify-spaces}

Following a suggestion of M.~Kreck, we use classification results
for highly connected manifolds
to identify some of the manifolds
from Sections~\ref{s-1-case}~and~\ref{z-2-case}
up to homeomorphism.
In all cases we have already computed the integral homology.
To apply the classification results, we need to know
the cup product structure in cohomology, as well as certain
characteristic classes.

\subsection{The torus case}

The spaces
$X_r$,~$Y_r$ and~$Z_r$ are $r$-connected and -- except for
$r=1$ -- it follows already from Poincar\'e duality and degree
reasons that $H^*(Z_r)$ is isomorphic as graded ring 
to the cohomology of a connected sum of products of spheres.

\subsubsection*{Case \texorpdfstring{$r=1$}{r=1}}

According to the construction given in Section~\ref{s-1-case},
$Z_1$ is the quotient of~$D^2\times T$ by the following
equivalence relation:
For points~$(x,g)$ with $x$~contained in the interior of the disk
no identification takes place.
The boundary~$S^1\subset D^2$ is divided into $4$~segments,
and for~$x$ in two opposite segments one identifies
one coordinate circle~$S^1\subset T=S^1\times S^1$
and for the other two segments the other coordinate circle.
For the $4$~points~$x$ separating the segments one identifies all of~$T$.
Replacing the disk by a square, one arrives at the usual
topological construction of the toric manifold~$S^2\times S^2$
where one looks at a $2$-sphere as a quotient of $[0,1]\times S^1$.
Hence,
\begin{equation}\label{Z1-spheres}
  Z_1\cong S^2\times S^2.
\end{equation}
Since we know $H_T^*(Z_1)$ to be free over~$R$,
we could alternatively
use Theorem~\ref{C-S-Z} to compute the cup product
in $H_T^*(Z_1)$ and $H^*(Z_1)=H_T^*(Z_1)\otimes_R\Z$
and then apply classification results for $4$-manifolds
(see~\cite{Freedman:82}) to conclude \eqref{Z1-spheres}.

\subsubsection*{Case \texorpdfstring{$r=2$}{r=2}}

We claim that the first Pontryagin class~$p_1(Z_2)$ vanishes.
Since the map $H^4(Z_2)\to H^4(Y_2)$ is injective, it suffices to show
that $p_1(Y_2)=0$. It follows from the fibration~$S^1\to Y_2\to X_2$
that the tangent bundle of~$Y_2$ is the Whitney sum of the
pull-back of the tangent bundle of~$X_2$ and the tangent bundle
along the fibres, and the latter plus a one-dimensional trivial
bundle is the pull-back of the vector bundle associated to the
above fibration (cf.~\cite{Szczarba:64}, for example).
The first Pontryagin class of this vector bundle vanishes
(since it is induced from the Hopf bundle over~$S^2$)
and so does $p_1(X_2)=0$ (since $X_2$ is an open subset of~$S^6$).
So by the product formula for Pontryagin classes (see~\cite{Milnor:74})
and the fact that $H^*(Y_2)$ has no $2$-torsion, one has $p_1(Y_2)=0$.
From~\cite{Wilkens:72} one gets
\begin{equation}
  Z_2\cong 3\connmult(S^3\times S^4).
\end{equation}

\subsubsection*{Cases \texorpdfstring{$r=4$}{r=4}~and~\texorpdfstring{$r=8$}{r=8}}

Here the cup product structure is clear as well (see above)
and one could prove the vanishing of certain characteristic
classes along the same line of arguments, but we do not know of a
classification result which allows to identify $Z_r$ in these
cases up to homeomorphism.

\subsection{The finite case}

The spaces~$X_r$,~$Y_r$ and~$Z_r$ are clearly $(r-1)$-connected
in this case.

\subsubsection*{Case \texorpdfstring{$r=1$}{r=1}}

It is elementary to see that
\begin{equation}
  Z_1\cong S^1\times S^1.
\end{equation}


For~$r\in\{2,4\}$ we consider the restriction of the 
$S^{r-1}$-action on~$Z_r$ to~$K=S^1\subset S^{r-1}$.
This action has $l=2^{r+1}$~fixed points. Since
$H^*(Z_r)$~and~$H^*(Z_r^K)$ are free $\Z$-modules
of the same rank, $H_K^*(Z_r)$ is free over~$H^*(BK)$.
As a consequence we get an injection
\begin{equation}
  H_K^*(Z_r)\to H_K^*\bigl(Z_r^K\bigr)\cong \Z[t]^{l}
\end{equation}
(which can be seen as a very special case of Theorem~\ref{C-S-Z}).

We write $b^{\langle n\rangle}$ as a shorthand for~$(b,\dots,b)\in\Z^n$.
As shown in~\cite[Sec.~2]{Puppe:01},
the cohomology $H^*(Z_r)$ can be described
as the graded algebra associated
to the following filtration of~$\Z^{l}$:
$\FF_0=\dots=\FF_{r-1}$ 
is the $\Z$-module generated by~$1^{\langle l\rangle}$,
$\FF_r=\dots=\FF_{2r-1}=\ker\sigma$ where
$\sigma\colon\Z^{l}\to\Z$ is given by~$\sigma (e_i)=(-1)^i$
for~$i=1,\dots,l$,
and $\FF_{2r}=\Z^{l}$.

The intersection form on~$H^r(Z_r)\cong\FF_r/\FF_{r-1}$ with
respect to the basis represented by
\begin{equation}
  v_i=\begin{cases}
    (0^{\langle i\rangle},1,1,0^{\langle l-i-2\rangle}) & \text{for odd~$i$,}\\
    (1^{\langle i\rangle},0^{\langle l-i\rangle}) & \text{for even~$i$,}\\
  \end{cases}
  \quad
  i=1,\dots,l-2,
\end{equation}
is given as a direct sum of~$(l-2)/2$~copies of the form~$\bigl[
  \begin{smallmatrix}
    0 & 1 \\
    1 & 0
  \end{smallmatrix}
\bigr]$.
Hence also for~$r\in\{2,4\}$
we have an isomorphism of graded rings
between the cohomology of~$Z_r$ and
that of a connected sum of $2^r-1$~copies of~$S^r\times S^r$.

\subsubsection*{Case \texorpdfstring{$r=2$}{r=2}}

It now follows from \cite{Freedman:82} that
\begin{equation}
  Z_2\cong 3\connmult(S^2\times S^2).
\end{equation}

\subsubsection*{Case \texorpdfstring{$r=4$}{r=4}}

As before, one can show that $p_1(Z_4)=0$ by using the
fibration~$S^3\to Y_4\to X_4$ and the facts that the pull-back of
the first Pontryagin class to the total space of the Hopf bundle
over $S^4$ vanishes (which can be seen e.g. from the Gysin
cohomology sequence of the sphere bundle) and that $p_1(X_4)=0$.
By \cite{Wall:62}~and~\cite{Schmitt:02} this together with the cup
product structure implies
\begin{equation}
  Z_4\cong 15\connmult(S^4\times S^4).
\end{equation}

\end{document}